\documentclass[a4paper, 11pt]{amsart}

\usepackage[ansinew]{inputenc}

\usepackage{etoolbox}
\makeatletter
\let\ams@starttoc\@starttoc
\makeatother
\usepackage[parfill]{parskip}
\makeatletter
\let\@starttoc\ams@starttoc
\patchcmd{\@starttoc}{\makeatletter}{\makeatletter\parskip\z@}{}{}
\makeatother

\usepackage{color,soul}
\definecolor{red}{rgb}{1,0,0}
\setulcolor{red}

\usepackage{amssymb}
\usepackage{amsmath}
\usepackage{amsthm}
\usepackage{fancyhdr}
\usepackage{esint}

\usepackage[colorlinks=true,linkcolor=blue,citecolor=blue,urlcolor=blue]{hyperref}
\usepackage{bookmark}

\swapnumbers
\newtheorem{lemma}{Lemma}[section]
\newtheorem{prop}[lemma]{Proposition}
\newtheorem{thm}[lemma]{Theorem}
\newtheorem{cor}[lemma]{Corollary}
\newtheorem{ass}[lemma]{Assumption}

\theoremstyle{definition}
\newtheorem{defn}[lemma]{Definition}
\newtheorem{example}[lemma]{Example}
\newtheorem{rem}[lemma]{Remark}

\renewcommand{\(}{\left(}
\renewcommand{\)}{\right)}
\renewcommand{\~}{\tilde}
\renewcommand{\-}{\bar}

\newcommand{\cn}{\colon}

\newcommand{\R}{\mathbb{R}}

\renewcommand{\S}{\mathbb{S}}

\renewcommand{\a}{\alpha}
\renewcommand{\b}{\beta}
\newcommand{\g}{\gamma}
\renewcommand{\d}{\delta}
\newcommand{\e}{\epsilon}
\renewcommand{\k}{\kappa}
\renewcommand{\l}{\lambda}

\newcommand{\vt}{\vartheta}
\renewcommand{\t}{\theta}
\newcommand{\s}{\sigma}
\newcommand{\p}{\varphi}

\newcommand{\T}{\Theta}

\newcommand{\del}{\partial}

\DeclareMathOperator{\graph}{graph}

\newcommand{\Def}{\begin{defn}}
\newcommand{\eDef}{\end{defn}}
\newcommand{\Prop}{\begin{prop}}
\newcommand{\eProp}{\end{prop}}
\newcommand{\Rem}{\begin{rem}}
\newcommand{\eRem}{\end{rem}}
\newcommand{\Lem}{\begin{lemma}}
\newcommand{\eLem}{\end{lemma}}
\newcommand{\Al}{\begin{align*}}
\newcommand{\eAl}{\end{align*}}
\newcommand{\eq}{\begin{equation}}
\newcommand{\eeq}{\end{equation}}
\newcommand{\ex}{\begin{example}}
\newcommand{\eex}{\end{example}}
\newcommand{\pf}{\begin{proof}}
\newcommand{\epf}{\end{proof}}
\newcommand{\Cor}{\begin{cor}}
\newcommand{\eCor}{\end{cor}}

\newcommand{\ra}{\rightarrow}
\newcommand{\hra}{\hookrightarrow}

\newcommand{\mc}{\mathcal}

\newcommand{\mrm}{\mathrm}

\newcommand{\hp}{\hphantom}

\begin{document}

\numberwithin{equation}{section}

\hyphenation{mani-folds mani-fold Rie-man-nian}

\title[The inverse mean curvature flow in warped cylinders]{The inverse mean curvature flow in warped cylinders of non-positive radial curvature}
\author{Julian Scheuer}
\subjclass[2000]{35J60, 53C21, 53C44, 58J05}
\keywords{curvature flows, inverse curvature flows, warped products}
\thanks{This work has been supported by the DFG}
\date{\today}
\address{Ruprecht-Karls-Universit\"at, Institut f\"ur Angewandte Mathematik, Im Neuenheimer Feld 294, 69120 Heidelberg, Germany}
\email{julian.scheuer@math.uni-freiburg.de}
\begin{abstract}
We consider the inverse mean curvature flow in smooth Riemannian manifolds of the form $([R_{0},\infty)\times\S^{n},\-{g})$ with metric
$\-{g}=dr^{2}+\vt^{2}(r)\s$ and non-positive radial sectional curvature. We prove, that for initial mean-convex graphs over $\S^{n}$ the flow exists for all times and remains a graph over $\S^{n}.$ Under weak further assumptions on the ambient manifold, we prove optimal decay of the gradient and that the flow leaves become umbilic exponentially fast. We prove optimal $C^{2}$-estimates in case that the ambient pinching improves.
\end{abstract}

\maketitle

\tableofcontents

\section{Introduction}

In this work we consider the inverse mean curvature flow (IMCF)
\eq \dot{x}=\frac{1}{H}\nu\eeq in an ambient manifold of the form $([R_{0},\infty)\times\S^{n},\-{g}).$ For the detailed set of assumptions see \ref{assA}, \ref{assB} and \ref{assC}. The motivation for this investigation arose from studying several works on the IMCF or on more general curvature flows in ambient spaces of constant or asymptotically constant sectional curvature, as Gerhardt's works in the Euclidean space,  \cite{Gerhardt:/1990} and \cite{Gerhardt:01/2014}, also compare \cite{Urbas:/1990}, or in the hyperbolic space, \cite{Gerhardt:11/2011}, for 1-homogeneous curvature functions and \cite{Scheuer:05/2015} for other homogeneities. In most of those works, long time existence of the flow as well as a smoothing effect, leading to roundness of the flow hypersurfaces, were established. For example, in the hyperbolic space the estimate for the second fundamental form \eq  |h^{i}_{j}-\d^{i}_{j}|\leq ce^{-\frac{2t}{n}}\eeq was proven in \cite{Scheuer:05/2015}. Depending on the specific setting, especially for curvature functions of higher homogeneity, additional assumptions like the convexity of the flow hypersurfaces had to be made. Similar convergence results for inverse curvature flows in the sphere were established in \cite{Gerhardt:/2015} and \cite{MakowskiScheuer:/2013}. Ding has also investigated the IMCF in the hyperbolic space and in rotationally symmetric spaces of Euclidean volume growth, \cite{Ding:01/2011}. For the sake of a complete presentation, we will not rule out this case in our proofs. Examples of ambient spaces of asymptotically constant sectional curvature can be found in \cite{BrendleHungWang:01/2016}. In many of those works, some of the a priori estimates have shown to be valid in more general settings, e.g. a bound on the gradient in finite time not assuming convexity, compare for example \cite[Lemma~3.6, (3.48)]{Gerhardt:01/2014} and a bound on the flow velocity in finite time as in \cite[Prop.~10]{BrendleHungWang:01/2016}. It is a natural question, whether the further estimates also hold in a more general setting, e.g. if one can prove long time existence and whether the flow hypersurfaces become umbilic as in case of constant sectional curvature spaces.

The IMCF has proven to be a very useful tool in various directions. In \cite{HuiskenIlmanen:/2001} Huisken and Ilmanen used a weak formulation of the IMCF to prove the Riemannian Penrose inequality in asymptotically flat manifolds. Other geometric inequalities have been proven also using more general curvature functions, e.g. Alexandrov-Fenchel inequalities as in  \cite{GeWangWu:04/2014}, also compare \cite{MakowskiScheuer:/2013}. Similar applications can be found in the work already mentioned above, \cite{BrendleHungWang:01/2016}, in which the IMCF was used to prove a Minkowski type inequality in the anti de Sitter Schwarzschild manifold. Of course we hope, and are optimistic, that the results of the present work allow applications in similar directions.

Now we describe in detail, what kind of ambient manifolds we allow in this work. The first assumption will be sufficient to prove long time existence of the flow.

\begin{ass}\label{assA}
Let $N=[R_{0},\infty)\times\S^{n},$ $n\geq 2,$ be equipped with its standard smooth structure and with the Riemannian metric 
\eq \label{RiemMetric} \-{g}=dr^{2}+\vt^{2}(r)\s_{\S^{n}},\eeq
where $\vt\in C^{\infty}([R_{0},\infty),\R_{+}),$ such that the radial sectional curvatures of $N$ are non-positive everywhere, in particular we assume
\eq \vt''\geq 0\ \wedge\ \vt'(r)>0\ \ \forall r>R_{0}.\eeq 
\end{ass}

For example, this assumption will be sufficient to prove long time existence of the IMCF in rotationally symmetric Hadamard manifolds $N$ with closed and mean-convex initial hypersurface, which is star-shaped around a point of rotational symmetry of $N.$ Then the only relevant part of the Hadamard manifold $N=(0,\infty)\times\S^{n},$ represented in geodesic polar coordinates, will be $[R_{0},\infty)\times\S^{n},$ where $R_{0}$ is small enough. Also note, that we do not impose a lower bound on the sectional curvature of $N.$

The following assumption will be sufficient to prove asymptotical roundness of the flow hypersurfaces. We will obtain an estimate of the form
\eq \label{umbilic} \left|h^{i}_{j}-\frac{\vt'}{\vt}\d^{i}_{j}\right|\leq \frac{ct}{\vt'}e^{-\frac{t}{n}}.\eeq

\begin{ass}\label{assB}
Suppose that for all $r_{0}>R_{0}$ there exist $c_{1},c_{2}>0,$ such that for all $r\geq r_{0}$
\eq \label{assB1} \vt''\vt\leq c_{1}\vt'^{2}\leq c_{2}(\vt''\vt+1)\eeq
and
\eq \label{assB2} \left|\frac{\vt'''}{\vt'}\right|\leq c\frac{\vt''}{\vt}.\eeq
\end{ass}

As we will describe in detail later, this assumption implies a bounded ratio of the sectional curvatures of $N$ in radial and in spatial direction. It says, that the space must not oscillate between hyperbolic and Euclidean behavior to quickly. However, and those are the remarkably new spaces compared to the works which already exist, the sectional curvatures of $N$ may oscillate between non-positive values infinitely often. Especially repeated transition from Euclidean to hyperbolic behavior is of interest.

The estimate (\ref{umbilic}) can be optimized (also compare Remark \ref{Remmainthm} for an explanation of 'optimality' in this context), if we impose improving pinching of $N,$ namely

\begin{ass}\label{assC}
Suppose that $N$ is asymptotically of constant curvature in the sense that there exist $c>0$ and $\l>2,$ such that
\eq \vt'\leq c\vt,\eeq
\eq |\t|\leq c\frac{\vt'^{2}}{\vt^{\l}}\eeq
and \eq |\t'|\leq c\frac{\vt'^{3}}{\vt^{1+\l}},\eeq
where \eq \t=\frac{\vt''}{\vt}+\frac{1-\vt'^{2}}{\vt^{2}}.\eeq
\end{ass}

In this case we are able to remove the extra $t$ in estimate (\ref{umbilic}). Beside the spaces of constant non-positive curvature, there are several other important manifolds, that satisfy the Assumptions \ref{assA}, \ref{assB} and \ref{assC}. The de Sitter-Schwarzschild manifolds with $\vt$ satisfying
\eq \vt'=\sqrt{1+\k\vt^{2}-m\vt^{1-n}},\ \k,m\geq 0,\eeq  are of this kind, cf. \cite[p.~3]{Brendle:06/2013}. We recapture the result by Brendle, cf. \cite{BrendleHungWang:01/2016} and also remove the extra factor $t^{2}$ in his estimate of the second fundamental form, cf. \cite[Prop.~16]{BrendleHungWang:01/2016}. In his work, Brendle already mentioned, that this would be possible, but he did not carry out the proof in this paper. The Reissner-Nordstrom manifolds, which satisfy
\eq \vt'=\sqrt{1-m\vt^{1-n}+q^{2}\vt^{2-2n}},\ m>0,\eeq
also satisfy our assumptions, cf. \cite[Ch. 5]{Brendle:06/2013}.  

Now we formulate the main theorem of this work.

\begin{thm} \label{mainthm}
(i) Let $(N,\-{g})$ satisfy Assumption \ref{assA}. Let $M\hra M_{0}\subset \mathring N$ be the embedding of a smooth, closed and mean convex hypersurface, which can be written as a graph over $\S^{n},$
\eq M_{0}=\graph u(0,\cdot),\eeq
where $u(t,x^{i})=r(0,\xi^{i}), (x^{i})\in \S^{n}, (\xi^{i})\in M.$
Then there is a unique smooth curvature flow
\eq x\cn [0,\infty)\times M\ra N\eeq
satisfying
\begin{align}\begin{split}\label{floweq} \dot{x}&=\frac{1}{H}\nu\equiv\frac{1}{\sum_{i=1}^{n}\k_{i}}\nu \\
								x(0)&=M_{0}, \end{split}\end{align}
where $\nu(t,\xi)$ is the outward normal to $M_{t}=x(t,M)$ at $x(t,\xi),$ $\k_{i}, 1\leq i\leq n,$ are the principal curvatures of $M_{t}$ in $x(t,\xi)$ and the leaves are graphs over $\S^{n},$
\eq M_{t}=\graph u(t,\cdot)\eeq
with a function
\eq u\in C^{\infty}([0,\infty)\times\S^{n},N).\eeq
(ii) Under Assumption \ref{assB}, the leaves $M_{t}$ of the flow become umbilical exponentially fast. There exists $c=c(N,M_{0}),$ such that
\eq \left|h^{i}_{j}-\frac{\vt'}{\vt}\d^{i}_{j}\right|\leq \frac{ct}{\vt'}e^{-\frac{t}{n}}\ \ \forall t\in[0,\infty).\eeq
Furthermore the gradient of the graph function $u$ satisfies the estimate
\eq \label{mainthm1} \|Du\|_{g}\leq \frac{c}{\vt'},\eeq where $g$ is the induced metric of the flow hypersurfaces, and we have
\eq \label{mainthm2} \sup\limits_{x\in\S^{n}}\vt'\|Du\|_{g}\ra 0,\eeq
if \eq\label{RescaledtoSphere} t\mapsto\frac{1}{\vt'^{2}}\(\vt^{-1}(\vt(R_{0})e^{\frac{t}{n}}\)\notin L^{1}([0,\infty)).\eeq
(iii) Under Assumption \ref{assC} we obtain the estimate
\eq \label{mainthm3}\left|h^{i}_{j}-\frac{\vt'}{\vt}\d^{i}_{j}\right|\leq \frac{c}{\vt'}e^{-\frac{t}{n}}\ \ \forall t\in[0,\infty).\eeq
\end{thm}

\Rem \label{Remmainthm}
Estimates (\ref{mainthm1}) and (\ref{mainthm3}) are exactly the ones, which are required to prove $C^{1}$ and $C^{2}$ boundedness of the rescaled surfaces 
\eq \~{u}=\vt(u)e^{-\frac{t}{n}}\eeq with respect to the metric of $\S^{n}.$ In general, this estimate can not be improved in this setting, i.e. (\ref{mainthm2}) does not hold in general, which can be seen by considering a geodesic sphere $\mc{S}$ in the hyperbolic space ($ \vt=\sinh$), written as a graph over another small geodesic sphere around a point $p_{0}$ not being the center of $\mc{S}.$ If (\ref{mainthm2}) would hold, then the gradient of the graph function arising from the evolution of $\mc{S}$ would decay to $0.$ But then
\eq e^{u-\frac{t}{n}}=2\vt e^{-\frac{t}{n}}+e^{-u-\frac{t}{n}}\eeq would converge to a constant, which is impossible, because the oscillation of the flow surfaces arising from $\mc{S}$ is positively constant with respect to $p_{0}$.
In the Euclidean case, however, (\ref{mainthm2}) holds, because the rescaling is
\eq \~{u}=u e^{-\frac{t}{n}}\eeq in this case,
where oscillations of the graphs are killed. From this behavior we see, that the choice of the sphere over which the surfaces are graphed is important to describe the flow behavior in an optimal way. In our general ambient manifold, however, we have no choice and thus can not expect a better decay than (\ref{mainthm1}). Also compare \cite{HungWang:09/2014} where it was shown that in general there is no such optimal center in the hyperbolic space and \cite{Scheuer:07/2016} for a proof that there is an optimal center for flows in $\R^{n+1}$ such that the oscillations of the hypersurfaces tend to zero without rescaling.
\eRem

\section{Setting, notation and general facts}
Now we state some general facts about hypersurfaces, especially those that can be written as graphs. We basically follow the description of \cite{Gerhardt:11/2011} and \cite{Scheuer:05/2015}, but restrict to Riemannian manifolds. For a detailed discussion we refer to \cite{Gerhardt:/2006}. \ \\
Let $N=N^{n+1}$ be Riemannian and $M=M^{n}\hookrightarrow N$ be a hypersurface. The geometric quantities of $N$ will be denoted by $(\-{g}_{\a\b}),$ $(\-{R}_{\a\b\g\d})$ etc., where greek indices range from $0$ to $n$. Coordinate systems in $N$ will be denoted by $(x^{\a}).$\ \\
Quantities for $M$ will be denoted by $(g_{ij}),$ $(h_{ij})$ etc., where latin indices range from $1$ to $n$ and coordinate systems will generally be denoted by $(\xi^{i}),$ unless stated otherwise.\ \\
Covariant differentiation will usually be denoted by indices, e.g. $u_{ij}$ for a function $u\colon M\rightarrow \R$, or, if ambiguities are possible, by a semicolon, e.g.
$h_{ij;k}.$ Usual partial derivatives will be denoted by a comma, e.g. $u_{i,j}.$\ \\
Let $x\colon M\hookrightarrow N$ be an embedding and $(h_{ij})$ be the second fundamental form, then we have the \textit{Gaussian formula}
\eq x^{\a}_{ij}=-h_{ij}\nu^{\a}, \eeq where $\nu$ is a differentiable normal, the \textit{Weingarten equation}
\eq \nu^{\a}_{i}=h^{k}_{i}x^{\a}_{k}, \eeq
the \textit{Codazzi equation} \eq h_{ij;k}-h_{ik;j}=\-{R}_{\a\b\g\d}\nu^{\a}x^{\b}_{i}x^{\g}_{j}x^{\d}_{k} \eeq
and the \textit{Gau\ss\ equation} \eq R_{ijkl}=(h_{ik}h_{jl}-h_{il}h_{jk})+\-{R}_{\a\b\g\d}x^{\a}_{i}x^{\b}_{j}x^{\g}_{k}x^{\d}_{l}. \eeq

Now assume that $N=(a,b)\times S_{0},$ where $S_{0}$ is compact Riemannian and that there is a Gaussian coordinate system $(x^{\a})$ such that
\eq d\-{s}^{2}=e^{2\psi}((dx^{0})^{2}+\-{g}_{ij}(x^{0},x)dx^{i}dx^{j}), \eeq
where $\-{g}_{ij}$ is a Riemannian metric, $x=(x^{i})$ are local coordinates for $S_{0}$ and $\psi\colon N\rightarrow \R$ is a function.\ \\
Let $M=\graph u_{|S_{0}}$ be a hypersurface
\eq M=\{(x^{0},x)\colon x^{0}=u(x), x\in S_{0}\}, \eeq
then the induced metric has the form \eq g_{ij}=e^{2\psi}(u_{i}u_{j}+\-{g}_{ij}) \eeq
with inverse \eq g^{ij}=e^{-2\psi}(\-{g}^{ij}-v^{-2}u^{i}u^{j}), \eeq
where $(\-{g}^{ij})=(\-{g}_{ij})^{-1},$ $u^{i}=\-{g}^{ij}u_{j}$ and \eq v^{2}=1+\-{g}^{ij}u_{i}u_{j}\equiv 1+|Du|^{2}.\eeq
We use, especially in the Gaussian formula, the normal 
\eq \label{outernormal} (\nu^{\a})=v^{-1}e^{-\psi}(1,-u^{i}). \eeq
Looking at $\a=0$ in the Gaussian formula, we obtain 
\eq e^{-\psi}v^{-1}h_{ij}=-u_{ij}-\-{\Gamma}^{0}_{00}u_{i}u_{j}-\-{\Gamma}^{0}_{0i}u_{j}-\-{\Gamma}^{0}_{0j}u_{i}-\-{\Gamma}^{0}_{ij} \eeq
and \eq e^{-\psi}\-{h}_{ij}=-\-{\Gamma}^{0}_{ij}, \eeq
where covariant derivatives are taken with respect to $g_{ij}.$

\subsection*{Rotationally symmetric spaces}

In this work we consider ambient spaces $N$ as described in Assumption \ref{assA}. We fix the coordinate system of $N,$ in which (\ref{RiemMetric}) holds, from now on. Then the Riemannian curvature tensor is given by 
\eq \label{RiemCurv}\bar{R}_{\a\b\g\d}=-\frac{\vt''}{\vt}\-{S}_{\a\b\g\d}+\t\-S_{\a'\b'\g'\d'}P^{\a'}_{\a}P^{\b'}_{\b}P^{\g'}_{\g}P^{\d'}_{\d},  \eeq
where \eq \-{S}_{\a\b\g\d}=\-{g}_{\a\g}\-{g}_{\b\d}-\-{g}_{\a\d}\-{g}_{\b\g},\eeq
\eq P^{\a'}_{\a}=\d^{\a'}_{\a}-r^{\a'}r_{\a}\eeq 
and \eq \t=\frac{\vt''}{\vt}+\frac{1-\vt'^{2}}{\vt^{2}},\eeq
compare the proof of \cite[Lemma~2.2]{BrendleHungWang:01/2016}.
Then, the Ricci tensor is given by
\eq \label{Ricci}\bar{R}_{\a\b}=-n\frac{\vt''}{\vt}\-{g}_{\a\b}+(n-1)\t\-{g}_{ij}.\eeq
 The second fundamental form of the coordinate slices $\{r=\mathrm{const}\}$ is given by
\eq \-{h}_{ij}=\frac{\vt'}{\vt}\-{g}_{ij}\eeq and thus
\eq \-{\k}_{i}=\frac{\vt'}{\vt}\ \ \forall 1\leq i\leq n,\eeq as the Riemannian curvature tensor reveals. 

If a hypersurface $M$ is a graph of a function $u$ in this coordinate system, we obtain
in this special situation
\eq h_{ij}v^{-1}=-u_{ij}+\-{h}_{ij},\eeq where covariant differentiation is performed with respect to the induced metric and
\eq \label{graphA} h_{j}^{i}=\frac{\vt'}{v\vt}\d^{i}_{j}+\frac{\vt'}{v^{3}\vt^{3}}u^{i}u_{j}-\frac{\~{g}^{ik}}{v\vt^{2}}u_{:kj}, \eeq
where $\~{g}^{ik}$ is the inverse of \eq \~{g}_{ik}=\p_{i}\p_{j}+\sigma_{ij},\eeq 
\eq \label{phi}\p(x)=\int_{R_{0}}^{u(x)}\vt^{-1}\eeq and covariant differentiation as well as index raising is performed with respect to the spherical metric $\sigma_{ij}$ (note the colon in $u_{:kj}$ to distinguish this derivative from the one with respect to the induced metric). For a proof of those formulas we refer to \cite[Thm. 3.13]{Scheuer:05/2015}.

\subsection*{Evolution equations}
The following is a standard observation on short time existence as it was already formulated in \cite[Remark 3.4]{Scheuer:05/2015}. 

\Rem \label{shorttime}
Looking at \cite[Thm. 2.5.19, Thm. 2.6.1]{Gerhardt:/2006}, under Assumption \ref{assA} we obtain short time existence of the flow on a maximal time interval $[0,T^{*}),$ $0<T^{*}\leq \infty.$ The flow exists at least as long as corresponding scalar flow of the graph function,
\eq \label{ScalarFlow} \frac{\del u}{\del t}=\frac{v}{H}, \eeq
where \eq u\cn[0,\-{T})\times\S^{n}\ra\R\eeq
and $\-{T}\leq T^{*},$
also compare \cite[Thm. 2.5.17]{Gerhardt:/2006} and \cite[p. 98-99]{Gerhardt:/2006}. Thus, it suffices to prove long time existence for (\ref{ScalarFlow}).
\eRem

Now we state the relevant evolution equations for the geometric quantities involved in the flow, as long as the hypersurfaces can be written as a graph over $\S^{n}.$ The flow velocity
\eq -\Phi=\frac{1}{H}\eeq
satisfies, compare \cite[Lemma 2.3.4]{Gerhardt:/2006} and (\ref{Ricci}),
\eq \label{EvPhi} \dot{\Phi}-\frac{1}{H^{2}}\Delta\Phi=\frac{\|A\|^{2}}{H^{2}}\Phi-\frac{n}{H^{2}}\frac{\vt''}{\vt}\Phi+\frac{(n-1)\t}{H^{2}}\|Du\|^{2}\Phi,\eeq
where $\dot{\Phi}=\frac{d}{dt}\Phi$ and $\Delta$ is the Laplace operator with respect to the induced metric of the flow hypersurfaces and where we also used $\|Du\|^{2}=|Du|^{2}v^{-2},$ compare \cite[equ.~(2.7.69)]{Gerhardt:/2006}.

\Lem
The second fundamental form satisfies the evolution equation
\begin{align}\begin{split}\label{EvA} \dot{h}^{i}_{j}-\frac{1}{H^{2}}\Delta h^{i}_{j}&=\frac{1}{H^{2}}g^{kl}\left(h_{kr}h^{r}_{l}-2\frac{\vt'}{\vt}h_{kl}+\frac{\vt'^{2}}{\vt^{2}}g_{kl}\right)h^{i}_{j}-\frac{1}{\vt^{2}}\frac{n}{H^{2}}h^{i}_{j}\\
						&\hp{=}-\frac{2}{H}\left(h^{i}_{k}-\frac{\vt'}{\vt}\d^{i}_{k}\right)h^{k}_{j}-\frac{n\t}{H^{2}}v^{-2}h^{i}_{j}+\frac{2\t}{H}v^{-2}\d^{i}_{j}\\
						&\hp{=}+\frac{\t}{H^{2}}\|Du\|^{2}h^{i}_{j}-\frac{2}{H^{3}}H_{j}H^{i}+\frac{n\t}{H^{2}}h^{mi}u_{m}u_{j}\\
						&\hp{=}+\frac{n\t}{H^{2}}h^{m}_{j}u_{m}u^{i}-\frac{2\t}{H^{2}}h^{mk}u_{m}u_{k}\d^{i}_{j}\\
						&\hp{=}+\frac{\t'}{H^{2}}v^{-1}(nu_{j}u^{i}-\|Du\|^{2}\d^{i}_{j})\\
						&\hp{=}+\frac{2\t}{H^{2}}v^{-1}\frac{\vt'}{\vt}\left(\|Du\|^{2}\d^{i}_{j}-nu^{i}u_{j}\right).\end{split}\end{align}
\eLem

\pf
According to \cite[Lemma 2.4.1]{Gerhardt:/2006} we have
\begin{align}\begin{split}\label{EvA1} \dot{h}^{i}_{j}-\frac{1}{H^{2}}\Delta h^{i}_{j}&=\frac{\|A\|^{2}}{H^{2}}h^{i}_{j}-\frac{2}{H}h_{kj}h^{ki}-\frac{2}{H^{3}}H_{j}H^{i}+\frac{1}{H^{2}}\-{R}_{\a\b}\nu^{\a}\nu^{\b}h^{i}_{j}\\
					&\hp{=}-\frac{1}{H^{2}}g^{kl}\-{R}_{\a\b\g\d}(x^{\a}_{m}x^{\b}_{k}x^{\g}_{r}x^{\d}_{l}h^{m}_{j}g^{ri}+x^{\a}_{m}x^{\b}_{k}x^{\g}_{j}x^{\d}_{l}h^{mi})\\
					&\hp{=}+\frac{2}{H^{2}}g^{kl}\-{R}_{\a\b\g\d}x^{\a}_{m}x^{\b}_{j}x^{\g}_{k}x^{\d}_{r}h^{m}_{l}g^{ri}\\
					&\hp{=}-\frac{2}{H}\-{R}_{\a\b\g\d}\nu^{\a}x^{\b}_{j}\nu^{\g}x^{\d}_{m}g^{mi}\\
					&\hp{=}+\frac{1}{H^{2}}g^{kl}\-{R}_{\a\b\g\d;\e}(\nu^{\a}x^{\b}_{k}x^{\g}_{l}x^{\d}_{j}x^{\e}_{m}+\nu^{\a}x^{\b}_{j}x^{\g}_{k}x^{\d}_{m}x^{\e}_{l})g^{mi}.\end{split}\end{align}					
There hold, using (\ref{RiemCurv}),
\eq \-{R}_{\a\b\g\d}x^{\a}_{m}x^{\b}_{j}x^{\g}_{k}x^{\d}_{r}=-\frac{\vt''}{\vt}(g_{mk}g_{jr}-g_{mr}g_{jk})+\t(\-{g}_{mk}\-{g}_{jr}-\-{g}_{mr}\-{g}_{jk}),\eeq

\eq \-{R}_{\a\b\g\d}\nu^{\a}x^{\b}_{j}\nu^{\g}x^{\d}_{m}=-\frac{\vt''}{\vt}g_{jm}+\t\|Du\|^{2}\-{g}_{jm}-\t v^{-2}u_{j}u_{m}\eeq
and

\begin{align}\begin{split} \-{R}_{\a\b\g\d;\e}\nu^{\a}x^{\b}_{k}x^{\g}_{l}x^{\d}_{j}x^{\e}_{m}&=-\t'v^{-1}(u_{l}u_{m}\-{g}_{kj}-u_{j}u_{m}\-{g}_{kl})-\t v^{-1}\frac{\vt'}{\vt}\-S_{mklj}\\
		&\hp{=}+\frac{\vt'}{\vt}v^{-1}(g_{mj}u_{l}u_{k}-g_{ml}u_{j}u_{k}\\
		&\qquad\qquad\quad+g_{kj}u_{m}u_{l}-g_{kl}u_{m}u_{j}).\end{split}\end{align}
We have \eq g_{ij}=u_{i}u_{j}+\-{g}_{ij},\eeq where one should note that we are using coordinates $(\xi^{i})\in M.$
Inserting those terms into (\ref{EvA1}) and rearranging terms, we obtain the result. We omit the tedious computation here.
\epf

We have the evolution equation for the graph function $u,$
\eq \label{Evu} \dot{u}-\frac{1}{H^{2}}\Delta u=\frac{2}{H}v^{-1}-\frac{\vt'}{\vt}\frac{n}{H^{2}}+\frac{\vt'}{\vt}\frac{1}{H^{2}}\|Du\|^{2},\eeq
cf. \cite[Lemma 3.3.2]{Gerhardt:/2006}, as well as the equation for $v=\sqrt{1+|Du|^{2}},$ compare \cite[Lemma 2.4.5]{Gerhardt:/2006} and \cite[Lemma 5.7]{Gerhardt:01/1996}. Note that in these two references there is a sign error, namely for example in \cite[Lemma~2.4.5, equ.~(2.4.28)]{Gerhardt:/2006} the term
\eq \dot{\Phi}F^{ij}\-R_{\a\b\g\d}\nu^{\a}x^{\b}_{i}x^{\g}_{j}x^{\d}_{k}r_{\e}x^{\e}_{m}g^{mk}v^{2}\eeq
should have a minus sign in front, or equivalently the indices $j$ and $k$ swapped. In our special situation we obtain

\begin{align}\begin{split}\label{Evv} \dot{v}-\frac{1}{H^{2}}\Delta v&=-\frac{1}{H^{2}}g^{kl}\left(h_{rk}h^{r}_{l}-2\frac{\vt'}{\vt}h_{kl}+\frac{\vt'^{2}}{\vt^{2}}g_{kl}\right)v-\frac{2}{H}\frac{\vt'}{\vt}(v-1)\\
							&\hp{=}-\frac{n\t}{H^{2}}\|Du\|^{2}v+\frac{1}{\vt^{2}H^{2}}\|Du\|^{2}v\\
							&\hp{=}-\frac{2}{H^{2}}v^{-1}v_{i}v^{i}+\frac{2}{H^{2}}\frac{\vt'}{\vt}v_{i}u^{i}.\end{split}\end{align}

\subsubsection*{Evolution equations with respect to the spherical metric}

In order to estimate the gradient, we need another description of (\ref{ScalarFlow}). We follow the method in \cite[Sec.~3]{Gerhardt:11/2011}.

\Rem
We obtain another formula for $h^{i}_{j}$ in terms of $\p$, compare \eqref{phi},
\eq h^{i}_{j}=v^{-1}\vt^{-1}(\vt'\d^{i}_{j}-(\s^{ik}-v^{-2}\p^{i}\p^{k})\p_{:kj}),\eeq
where differentiation and index raising are performed with respect to $\s_{ij},$ cf. \cite[(3.26)]{Gerhardt:11/2011}. We obtain
\eq \label{phiScalarFlow} \frac{\del}{\del t}\p=\frac{v}{\vt H}=\cn\frac{v}{\~{H}}.\eeq
There  holds
\eq g_{ij}=u_{i}u_{j}+\vt^{2}\s_{ij}=\vt^{2}(\p_{i}\p_{j}+\s_{ij})\equiv\vt^{2}\~{g}_{ij}.\eeq
Note that
\eq |Du|^{2}=\vt^{-2}\s^{ij}u_{i}u_{j}=\s^{ij}\p_{i}\p_{j}\equiv|D\p|^{2},\eeq
as well as
\eq \label{phigraphA} \~{h}^{l}_{k}=-v^{-1}\~{g}^{lr}\p_{:rk}+v^{-1}\vt'\d^{l}_{k},\eeq
where $(\~{g}^{lr})=(\~{g}_{lr})^{-1}.$
\eRem

\Lem\label{phiDer}
The various quantities and tensors in (\ref{phiScalarFlow}) satisfy
\eq v_{i}=v^{-1}\p_{:ki}\p^{k},\eeq
\eq {\~{g}^{lr}}_{:i}=2v^{-3}v_{i}\p^{l}\p^{r}-v^{-2}({\p_{:}}^{l}_{i}\p^{r}+\p^{l}{\p_{:}}^{r}_{i})\eeq
and
\eq \~{h}^{l}_{k:i}=-\frac{v_{i}}{v}\~{h}^{l}_{k}-v^{-1}({\~{g}^{lr}}_{:i}\p_{:rk}+\~{g}^{lr}\p_{:rki}-\vt''\vt\p_{i}\d^{l}_{k}),\eeq
where all indices refer to $\s_{ij}.$
\eLem

\pf
This is a straightforward computation in any of the cases. Just have in mind, that $\vt=\vt(u),$ such that $\vt_{i}=\vt'u_{i}=\vt'\vt\p_{i}.$
\epf

We also need a version of the parabolic equation solved by $u$ with respect to the spherical metric.

\Lem \label{EvuB}
Let $u_{:ij}$ denote the second derivatives of $u$ with respect to $\s_{ij},$ $u_{;ij}$ those with respect to $\-{g}_{ij}=\vt^{2}\s_{ij}$ and $u_{ij}$ those with respect to $g_{ij}.$ Then there hold
\eq u_{:ij}=v^{2}u_{ij}+\frac{\vt'}{\vt}(2u_{i}u_{j}-\vt^{2}|Du|^{2}\s_{ij})\eeq
and
\eq \label{EvuB1}\mc{L}u:=\frac{\del}{\del t}u-\frac{1}{\~{H}^{2}}\~{g}^{kr}u_{:kr}=\frac{2v}{H}-\frac{n}{H^{2}}\frac{\vt'}{\vt}-\frac{1}{H^{2}}\frac{\vt'}{\vt}v^{-2}\s^{ij}\p_{i}\p_{j}.\eeq
\eLem

\pf
By \cite[Lemma 2.7.6]{Gerhardt:/2006} there holds
\begin{align}\begin{split} u_{ij}&=v^{-2}u_{;ij}=v^{-2}\(u_{i,j}-\frac{1}{2}\-{g}^{lm}(\-{g}_{il,j}+\-{g}_{jl,i}-\-{g}_{ij,l})u_{m}\)\\
						&=v^{-2}\Big(u_{i,j}-\frac{1}{2}\s^{lm}(\s_{il,j}+\s_{jl,i}-\s_{ij,l})u_{m}\\
						&\qquad\qquad-\frac{\vt'}{\vt}u^{l}(u_{j}\s_{il}+u_{i}\s_{jl}-u_{l}\s_{ij})\Big)\\
						&=v^{-2}u_{:ij}-2v^{-2}\frac{\vt'}{\vt}u_{i}u_{j}+v^{-2}\vt'\vt|Du|^{2}\s_{ij}.\end{split}\end{align}
Thus we have
\begin{align}\begin{split} \mc{L}u&=\frac{v}{H}-\frac{1}{\~{H}^{2}}\~{g}^{kr}u_{:kr}\\
						&=\frac{v}{H}-\frac{1}{\~{H}^{2}}\~{g}^{kr}\left(-vh_{kr}+v^{2}\-{h}_{kr}+2\frac{\vt'}{\vt}u_{k}u_{r}-\vt'\vt|Du|^{2}\s_{kr}\right)\\
						&=2\frac{v}{H}-\frac{nv^{2}}{\~{H}^{2}}\vt'\vt+\frac{v^{2}}{\~{H}^{2}}\frac{\vt'}{\vt}\~{g}^{kr}u_{k}u_{r}-\frac{2}{\~{H}^{2}}\frac{\vt'}{\vt}\~{g}^{kr}u_{k}u_{r}\\
						&\hp{=}+\frac{n}{\~{H}^{2}}\vt'\vt|Du|^{2}-\frac{1}{\~{H}^{2}}\frac{\vt'}{\vt}\~{g}^{kr}u_{k}u_{r}|Du|^{2}\\
						&=2\frac{v}{H}-\frac{n}{H^{2}}\frac{\vt'}{\vt}-\frac{1}{H^{2}}\frac{\vt'}{\vt}v^{-2}|D\p|^{2}.\end{split}\end{align}
\epf

Now we state the evolution equation for the gradient, for a proof we refer to \cite[Prop.~3.3]{BrendleHungWang:01/2016} and remember (\ref{phigraphA}).
\Lem\label{EvGrad}
Define
\eq F=\frac{\vt H}{v}.\eeq
Then the function
\eq w=\frac{1}{2}|D\p|^{2}\eeq
satisfies
\begin{align}\begin{split}\label{EvGrad1} \mc{L}w&=-\frac{1}{F^{2}}\frac{\del F}{\del\p_{i}}w_{i}-\frac{2(n-1)}{v^{2}F^{2}}w-\frac{2n}{v^{2}F^{2}}\vt''\vt w-\frac{1}{v^{2}F^{2}}\~{g}^{ij}{\p_{:}}^{k}_{i}\p_{:kj}\\
			&=-\frac{2n}{\~{H}^{2}}\vt''\vt w-\frac{2(n-1)}{\~{H}^{2}}w-\frac{1}{\~{H}^{2}}\~{g}^{ij}{\p_{:}}^{k}_{i}\p_{:kj}\\
			&\hp{=}+\frac{2}{\~{H}}v^{-1}w_{i}\p^{i}+\frac{2}{\~{H}^{2}}v^{-4}(\p^{i}w_{i})^{2}-\frac{2}{\~{H}^{2}}v^{-2}w_{i}w^{i}.\end{split}\end{align}
\eLem

\pagebreak

\section{Long time existence}

\subsection*{Barriers}

\Prop
Under Assumption \ref{assA} consider (\ref{floweq}) with $x(0)=S_{r_{0}}=\{x^{0}=r_{0}\},$ $r_{0}>R_{0}.$ Then the corresponding flow $x=x(t,\xi)$ exists for all times. The leaves $M_{t}=x(t,M)$ are constant graphs and
\eq x^{0}(t,M)=\Theta(t,r_{0}),\eeq
where $\T$ solves the initial value problem
\begin{align}\begin{split} \frac{d}{dt}\T&=\frac{\vt}{n\vt'}(\T)\\
				\T(0,r_{0})&=r_{0}.\end{split}\end{align}
Furthermore
\eq \lim\limits_{t\ra\infty}\T=\infty.\eeq
\eProp

\pf
In view of $\vt'\geq \vt'(r_{0}),$ $\vt$ is injective on $[r_{0},\infty).$ Define 
\eq \T(t)=\vt^{-1}\left(\vt(r_{0})e^{\frac{t}{n}}\right), 0\leq t<\infty.\eeq
This function has the desired properties.
\epf

The following lemma is standard, compare \cite[Lemma 6]{BrendleHungWang:01/2016}.

\Lem \label{BarriersA}
The solution $u$ of (\ref{ScalarFlow}) satisfies
\eq \T(t,\inf u(0,\cdot))\leq u(t,\cdot)\leq \T(t,\sup u(0,\cdot)),\eeq
as long as the flow hypersurfaces are graphs over $\S^{n}.$
\eLem

\Lem\label{BarriersB}
Under Assumption \ref{assA} let $\T_{i}(t)=\vt^{-1}\left(c_{i}e^{\frac{t}{n}}\right)$ be the spherical solutions of (\ref{ScalarFlow}) with $c_{2}> c_{1}>\vt(R_{0}).$ Then for all $t\in[0,\infty)$ there exists $\a\in(c_{1},c_{2}),$ such that
\eq 1\leq \frac{\vt'(\T_{2}(t))}{\vt'(\T_{1}(t))}\leq \frac{c_{2}-c_{1}}{c_{1}}\frac{\vt''\vt}{\vt'^{2}}\left(\vt^{-1}\left(\a e^{\frac{t}{n}}\right)\right).\eeq
\eLem

\pf
Let $t\in[0,\infty)$ be fixed. Then the mean value theorem implies the existence of $\a\in(c_{1},c_{2}),$ such that
\begin{align}\begin{split} &\log\vt'\left(\vt^{-1}\left(c_{2}e^{\frac{t}{n}}\right)\right)-\log\vt'\left(\vt^{-1}\left(c_{1}e^{\frac{t}{n}}\right)\right)\\
					=&\frac{\vt''}{\vt'^{2}}\left(\vt^{-1}\left(\a e^{\frac{t}{n}}\right)\right)e^{\frac{t}{n}}(c_{2}-c_{1})=\frac{\vt''\vt}{\vt'^{2}}\left(\vt^{-1}\left(\a e^{\frac{t}{n}}\right)\right)\frac{c_{2}-c_{1}}{\a}.\end{split}\end{align}
\epf

\subsection*{Gradient and curvature estimates}

\Lem\label{GradBound}
Under Assumption \ref{assA} let $u$ be the solution of (\ref{ScalarFlow}). Then the quantity
\eq \sup\limits_{x\in\S^{n}}v(\cdot,x)\eeq
is non-increasing.
\eLem

\pf
This follows from the maximum principle applied to (\ref{EvGrad1}).
\epf

In order to estimate the curvature function $H$ from below, we follow the method in \cite[Prop.~3.4]{BrendleHungWang:01/2016}.

\Lem\label{HBound}
Under Assumption \ref{assA} let $u$ be the solution of (\ref{ScalarFlow}). Then the function $\frac{\del}{\del t}\p$ is bounded from above. In particular we have
\eq H\geq c\frac{v}{\vt}.\eeq
\eLem

\pf
Define, as above, $F=\frac{\~{H}}{v}=\frac{\vt H}{v}.$ Then
\eq z:=\frac{\del}{\del t}\p=\frac{1}{F}.\eeq
Thus
\begin{align}\begin{split} \frac{\del}{\del t}z&=-\frac{1}{F^{2}}\left(\frac{\del F}{\del \p_{ij}}z_{:ij}+\frac{\del F}{\del \p_{i}}z_{i}+\frac{\del F}{\del \p}z\right )\\
						&=-\frac{1}{F^{2}}\left(v^{-1}\~{g}^{kl}\frac{\del\~{h}_{kl}}{\del\p_{ij}}z_{:ij}+\frac{\del F}{\del\p_{i}}z_{i}+v^{-1}\~{g}^{kl}\frac{\del\~{h}_{kl}}{\del\p}z\right)\\
						&=\frac{1}{F^{2}}v^{-2}\~{g}^{ij}z_{:ij}-\frac{1}{F^{2}}\frac{\del F}{\del\p_{i}}z_{i}-\frac{n}{F^{2}}v^{-2}\vt''\vt z.\end{split}\end{align}
The maximum principle implies the claim.
\epf

The last step in proving long time existence is to bound the principal curvatures.

\Lem\label{KappaBound}
Under Assumption \ref{assA} let $u$ be a solution of (\ref{ScalarFlow}). Then there exists $c=c(\-{T},N,M_{0}),$ such that the principal curvatures of the flow hypersurfaces $\k_{i}, 1\leq i\leq n,$ satisfy the pointwise estimate
\eq \k_{i}\leq c\frac{\vt'}{\vt}.\eeq
If Assumption \ref{assB} is satisfied additionally, then $c$ does not depend on $\-{T}.$
\eLem

\pf
Define
\eq \zeta=\sup\{h_{ij}\xi^{i}\xi^{j}\cn\|\xi\|=1\},\eeq
\eq w=\zeta\frac{\vt}{\vt'}\eeq
and \eq z=\log\zeta+\log\frac{\vt}{\vt'}+f(v),\eeq
where $f$ is yet to be determined. We want to bound $z.$ Thus suppose
\eq \sup\limits_{[0,T]\times M}z=z(t_{0},\xi_{0}), T<T^{*}, t_{0}>0.\eeq
In $(t_{0},\xi_{0})$ introduce Riemannian normal coordinates around $\xi_{0},$ such that
\eq g_{ij}=\d_{ij}, h_{ij}=\k_{i}\d_{ij}, \k_{1}\leq\dots\leq \k_{n}.\eeq
Without loss of generality we may assume, that in $(t_{0},\xi_{0})$ there holds
\eq h^{n}_{n}=\k_{n},\eeq
cf. the proof of Lemma 4.4 in \cite{Gerhardt:11/2011}. Using (\ref{EvA}), (\ref{Evu}) and (\ref{Evv}), we obtain

\begin{align}\begin{split}\label{KappaBound1} \dot{z}-\frac{1}{H^{2}}\Delta z&\leq\frac{1}{H^{2}}g^{kl}\left(h_{kr}h^{r}_{l}-2\frac{\vt'}{\vt}h_{kl}+\frac{\vt'^{2}}{\vt^{2}}g_{kl}\right)(1-f'v)-\frac{n}{\vt^{2}H^{2}}\\
				&\hp{=}-\frac{2}{H}\left(\k_{n}-\frac{\vt'}{\vt}\right)-\frac{n\t}{H^{2}}v^{-2}+\frac{2\t}{H}v^{-2}\k_{n}^{-1}+\frac{c|\t|}{H^{2}}\|Du\|^{2}\\
				&\hp{=}+\frac{c|\t'|}{H^{2}}\k_{n}^{-1}\|Du\|^{2}-\frac{2}{H^{3}}H_{n}H^{n}\k_{n}^{-1}+\frac{1}{H^{2}}\|D(\log h^{n}_{n})\|^{2}\\
				&\hp{=}-\frac{2f'}{H}\frac{\vt'}{\vt}(v-1)+cf'\frac{|\t|}{H^{2}}\|Du\|^{2}+\frac{f'}{\vt^{2}H^{2}}\|Du\|^{2}v\\
				&\hp{=}-\frac{2f'}{H^{2}}v^{-1}v_{i}v^{i}+\frac{2f'}{H^{2}}\frac{\vt'}{\vt}v_{i}u^{i}-\frac{f''}{H^{2}}v_{i}v^{i}\\
				&\hp{=}+\left(\frac{\vt'}{\vt}-\frac{\vt''}{\vt'}\right)\frac{2}{H}v^{-1}-\left(\frac{\vt'}{\vt}-\frac{\vt''}{\vt'}\right)\frac{\vt'}{\vt}\frac{n}{H^{2}}\\
				&\hp{=}+\left(\frac{\vt'^{2}}{\vt^{2}}-\frac{\vt''}{\vt}\right)\frac{1}{H^{2}}\|Du\|^{2}-\left(\frac{\vt'}{\vt}-\frac{\vt''}{\vt'}\right)'\frac{1}{H^{2}}\|Du\|^{2}.\end{split}\end{align}
Now choose $0<\a<v^{-1}$ and 
\eq f(v)=\log\frac{1}{\frac{1}{v}-\a}.\eeq
Using
\eq |\t|\leq \frac{\vt''}{\vt}+c\frac{\vt'^{2}}{\vt^{2}}=\frac{\vt'}{\vt}\left(\frac{\vt''}{\vt'}+c\frac{\vt'}{\vt}\right)\leq c(\-{T},N,M_{0})\frac{\vt'^{2}}{\vt^{2}}\eeq
and \eq |\t'|\leq c(\-{T},N,M_{0})\frac{\vt'^{3}}{\vt^{3}},\eeq 
where $c=c(N,M_{0})$ in case of Assumption \ref{assB}, 
\eq v_{i}=-v^{2}h^{k}_{i}u_{k}+v\frac{\vt'}{\vt}u_{i},\eeq
cf. \cite[(5.29)]{Gerhardt:01/1996},
as well as
\eq (\log h^{n}_{n})_{k}=-\left(\frac{\vt'}{\vt}-\frac{\vt''}{\vt'}\right)u_{i}-f'v_{i}\eeq
at $(t_{0},\xi_{0}),$ we find

\begin{align}\begin{split} 0&\leq -\frac{1}{H^{2}}g^{kl}\left(h_{kr}h^{r}_{l}-2\frac{\vt'}{\vt}h_{kl}+\frac{\vt'^{2}}{\vt^{2}}g_{kl}\right)\frac{\a v}{1-\a v}-\frac{2}{H}\left(\k_{n}-\frac{\vt'}{\vt}\right)\\
				&\hp{=}+\frac{c}{H^{2}}\frac{\vt'^{2}}{\vt^{2}}(1+w^{-1}+w)+\frac{1}{H^{2}}(f'^{2}-2f'v^{-1}-f'')v_{i}v^{i}\\
				&\leq -\frac{1}{H^{2}}\frac{\vt'^{2}}{\vt^{2}}\left(\frac{\a v}{1-\a v}(w-1)^{2}+c+cw+cw^{-1}\right),\end{split}\end{align}
from which we conclude
\eq w(t_{0},\xi_{0})\leq c=c(\-{T},N,M_{0}),\eeq where in case of Assumption \ref{assB} $c$ does not depend on $\-{T}.$
Thus we have proven the lemma.
\epf

Now we prove long time existence of the flow in case of Assumption \ref{assA}.

\begin{thm}
Let $N=[R_{0},\infty)\times \S^{n}$ be equipped with the warped product metric
\eq \-{g}=dr^{2}+\vt^{2}(r)\s,\eeq
such that $\vt\in C^{\infty}([R_{0},\infty),\R_{+})$ and such that \eq\vt''\geq 0\ \wedge\ \vt'_{|\mathring N}>0.\eeq Let $M_{0}\subset \mathring{N}$ be a smooth, closed and mean-convex embedded hypersurface, written as a graph over $\S^{n}.$ Then the inverse mean curvature flow
\begin{align}\begin{split} \dot{x}&=\frac{1}{H}\nu\\
						x(0)&=M_{0}\end{split}\end{align}
exists for all times.
\end{thm}

\pf
The corresponding scalar solution $u$ of (\ref{ScalarFlow}) satisfies the parabolic equation
\eq \frac{\del u}{\del t}=\frac{v}{H}=\frac{v}{n\frac{\vt'}{v\vt}+\frac{\vt'}{v^{3}\vt^{3}}u^{i}u_{i}-\frac{1}{v\vt^{2}}\~{g}^{ik}u_{:ik}}\equiv G(x,u,Du,D^{2}u),\eeq 
where indices are understood with respect to $\s_{ij},$ cf. (\ref{graphA}). Using the a priori estimates of this chapter, we see that $G$ is, during finite time, a uniformly parabolic operator, which satisfies the requirements to apply the regularity results of Krylov-Safonov, \cite[Par. 5.5]{Krylov:/1987}, to obtain
\eq |u|_{2,\a,\S^{n}}\leq c(\-{T},N,M_{0})\eeq
and \cite[Thm. 2.5.9]{Gerhardt:/2006} for higher order estimates. By compactness, $u$ may be extended beyond any finite $\-{T},$ also compare the proof of \cite[Lemma 2.6.1]{Gerhardt:/2006}. The result follows from the explanations in Remark \ref{shorttime}.
\epf

\section{Long time behavior}

In this section, we successively prove and improve the decay estimates, first for the gradient.

\Lem \label{GradDecay}
Under the assumptions \ref{assA} and \ref{assB} let $u$ be the solution of (\ref{ScalarFlow}). Then there exist $\g>0$ and $c=c(\g,N,M_{0}),$ such that
\eq |D\p|^{2}\leq \frac{c}{\vt'^{\g}}.\eeq
\eLem

\pf
Define \eq w=\frac{1}{2}|D\p|^{2}\eeq
 and \eq z=f(u)w,\eeq
where $f>0$ will be determined later with $f'\geq0.$ Suppose for $0<T<\infty,$ that
\eq \sup\limits_{[0,T]\times\S^{n}}z=z(t_{0},x_{0})\geq 1,\quad t_{0}>0.\eeq
Then from (\ref{EvuB1}) and (\ref{EvGrad1}) we obtain at $(t_{0},x_{0}),$ also using
\eq w_{i}=-\frac{f'}{f}w\vt\p_{i},\eeq

\begin{align}\begin{split}\label{GradDecay1}	\mc{L}z&=f\mc{L}w+z\frac{f'}{f}\mc{L}u-\frac{f''}{\~{H}^{2}}\~{g}^{kr}u_{k}u_{r}w-\frac{2}{\~{H}^{2}}\~{g}^{kr}w_{k}f_{r} \\
					&=-\frac{2n}{\~{H}^{2}}\vt''\vt z-\frac{2(n-1)}{\~{H}^{2}}z-\frac{f}{\~{H}^{2}}\~{g}^{ij}{\p_{:}}^{k}_{i}\p_{:kj}-\frac{4f'}{\~{H}}v^{-1}\vt w^{2}\\
					&\hp{=}+\frac{8}{\~{H}^{2}}v^{-4}\frac{f'^{2}}{f^{2}}\vt^{2}w^{3}z-\frac{4}{\~{H}^{2}}v^{-2}\frac{f'^{2}}{f^{2}}\vt^{2}w^{2}z+\frac{f'}{f}\frac{2v}{\~{H}^{2}}\vt\~{H}z\\
					&\hp{=}-\frac{f'}{f}\frac{n}{\~{H}^{2}}\vt\vt'z-\frac{f'}{f}\frac{1}{\~{H}^{2}}\vt\vt'\|Du\|^{2}z\\
					&\hp{=}+\frac{2}{\~{H}^{2}}\frac{f'^{2}}{f^{2}}\vt^{2}\|Du\|^{2}z-\frac{f''}{f}\frac{1}{\~{H}^{2}}\vt^{2}\|Du\|^{2}z. \end{split}\end{align}
We have
\eq \~{H}=\frac{n}{v}\vt'-v^{-1}\~{g}^{ij}\p_{:ij}.\eeq
Substituting this into (\ref{GradDecay1}), we obtain at $(t_{0},x_{0})$

\begin{align}\begin{split}\label{GradDecay2} 0&\leq \frac{2n}{\~{H}^{2}}z\left(-\vt''\vt-\frac{n-1}{n}+\frac{f'}{2f}\vt\vt'\right)-\frac{1}{\~{H}^{2}}\~{g}^{ij}{\p_{:}}^{k}_{i}\p_{:kj}f\\
				&\hp{=}-\frac{f'}{f}\frac{2}{\~{H}^{2}}\vt\~{g}^{ij}\p_{:ij}z-\frac{4}{\~{H}}v^{-1}\frac{f'}{f}\vt wz+\frac{8}{\~{H}^{2}}v^{-4}\frac{f'^{2}}{f^{2}}\vt^{2}w^{3}z\\
				&\hp{=}-\frac{2}{\~{H}^{2}}\frac{f'}{f}v^{-2}\vt\vt' wz+\frac{4}{\~{H}^{2}}\frac{f'^{2}}{f^{2}}\vt^{2}v^{-2}wz-\frac{2}{\~{H}^{2}}\frac{f''}{f}\vt^{2}v^{-2}wz.\end{split}\end{align}
Now we apply
\eq ab\leq \frac{\e}{2}a^{2}+\frac{1}{2\e}b^{2},\quad \e>0,\eeq
with the values
\eq a=|\~{g}^{ij}\p_{:ij}|,\quad b=1,\quad \e=\frac{\d}{w\vt'},\quad \d>0\eeq
and choose
\eq f(u)=\vt'^{\g}, \g>0.\eeq
Then (\ref{GradDecay2}) becomes

\begin{align}\begin{split}\label{GradDecay3} 0&\leq\mc{L}z\leq \frac{2n}{\~{H}^{2}}z\left(\left(\frac{\g}{2}-1\right)\vt''\vt-\frac{n-1}{n}+c\g^{2}v^{-4}(w^{3}+w)\vt''\vt\right)\\
			&\hp{=}-\frac{1}{\~{H}^{2}}\~{g}^{ij}{\p_{:}}^{k}_{i}\p_{:kj}\vt'^{\g}+\frac{\g\d}{\~{H}^{2}}\frac{\vt''\vt}{\vt'^{2}}|\~{g}^{ij}\p_{:ij}|^{2}\vt'^{\g}+\frac{\g}{\d\~{H}^{2}}\vt''\vt wz.\end{split}\end{align}
Using
\eq \~{g}^{ij}{\p_{:}}^{k}_{i}\p_{:kj}\geq \frac{1}{n}|\~{g}^{ij}\p_{:ij}|^{2},\eeq
we obtain a contradiction for small $\g>0$ and $\d=1.$
\epf

\Lem\label{HDecay}
Under the assumptions \ref{assA} and \ref{assB} let $u$ be the solution of (\ref{ScalarFlow}). Then there exists $c=c(N,M_{0}),$ such that
\eq H\frac{\vt}{\vt'}\leq c\eeq
and 
\eq \limsup\limits_{t\ra\infty}\sup_{M}H\frac{\vt}{\vt'}\leq n,\eeq
if $\vt'$ is unbounded.
\eLem

\pf
Define
\eq w=H\frac{\vt}{\vt'}=\cn Hf(u).\eeq
We have \eq f'=1-\frac{\vt''}{\vt'}f\eeq
and \eq f''=2\frac{\vt''^{2}}{\vt'^{2}}f-\frac{\vt'''}{\vt'}f-\frac{\vt''}{\vt'}.\eeq
Using (\ref{EvPhi}) and (\ref{Evu}), we obtain

\begin{align}\begin{split} \dot{w}-\frac{1}{H^{2}}\Delta w&=-\frac{\|A\|^{2}}{H^{2}}w+\frac{n}{H^{2}}\frac{\vt''}{\vt}w-\frac{(n-1)\t}{H^{2}}\|Du\|^{2}w\\
			&\hp{=}-\frac{2}{H^{3}}H_{i}H^{i}f+\frac{2f'}{v}-\frac{nf'}{H}\frac{\vt'}{\vt}+\frac{f'}{H}\frac{\vt'}{\vt}\|Du\|^{2}\\
			&\hp{=}-\frac{f''}{H}\|Du\|^{2}-\frac{2}{H^{2}}H_{i}f^{i}\\
			&=-\frac{\|A\|^{2}}{H^{2}}w+\frac{2}{v}+\frac{n}{w}\|Du\|^{2}-\frac{n}{w}-\frac{n-1}{\vt^{2}H^{2}}\|Du\|^{2}w\\
			&\hp{=}+\frac{\vt''}{\vt'}\frac{1}{H}\left(2n-\frac{2}{v}w-(n-1)\|Du\|^{2}\right)\\
			&\hp{=}-\frac{\|Du\|^{2}}{H^{2}}\left(2\frac{\vt''^{2}}{\vt'^{2}}-\frac{\vt'''}{\vt'}\right)w-\frac{2}{H^{3}}H_{i}H^{i}f-\frac{2}{H^{2}}H_{i}f^{i}\\
			&\leq -\left(\sqrt{\frac{w}{n}}-\sqrt{\frac{n}{w}}\right)^{2}+\|Du\|^{2}\left(\frac{n}{w}-2\frac{v}{v+1}-\frac{n-1}{w\vt'^{2}}\right)\\
			&\hp{=}+\frac{\vt''}{\vt'}\frac{1}{H}\|Du\|^{2}\left(2n\frac{v}{v+1}-(n-1)+c+c\frac{\vt''\vt}{\vt'^{2}}\right)\\
			&\hp{=}-\frac{\vt''}{\vt'}\frac{2}{H}v^{-1}(w-n)-\frac{2}{H^{3}}H_{i}H^{i}f-\frac{2}{H^{2}}H_{i}f^{i}.\end{split}\end{align}
Thus the function
\eq \~{w}=\sup\limits_{\xi\in M}(w(\cdot,\xi)-n)\eeq
is uniformly bounded by a constant $c=c(N,M_{0}),$ using assumptions \ref{assA} and \ref{assB}. If $\vt'\ra\infty,$ then by Lemma \ref{GradDecay} we have $\|Du\|^{2}\ra 0,$ which shows that $\~{w}$ is eventually strictly decreasing on any of the sets $\{\~{w}\geq \e>0\}.$ \cite[Lemma 4.2]{Scheuer:05/2015} implies
\eq \limsup\limits_{t\ra\infty}\~{w}\leq 0.\eeq
\epf

\Cor\label{GradDecayB}
Under the assumptions \ref{assA} and \ref{assB} let $u$ be the solution of (\ref{ScalarFlow}).Then, if $\vt'$ is uniformly bounded, there exist $c,\mu>0$ depending on $N$ and $M_{0},$ such that
\eq |D\p|^{2}\leq ce^{-\mu t}\ \ \forall t\in[0,\infty).\eeq
\eCor

\pf
We know at this moment, that
\eq \frac{1}{\vt H}\geq \frac{c}{\vt'}.\eeq
Thus the result follows from (\ref{EvGrad1}) immediately.
\epf

We finish our preparatory decay results by treating the mean curvature from below.

\Lem \label{HGrowth}
Under the assumptions \ref{assA} and \ref{assB} let $u$ be the solution of (\ref{ScalarFlow}).Then there exists $c=c(N,M_{0}),$ such that
\eq \frac{\vt'}{\vt H}\leq c\eeq
and in case that $\vt'$ is unbounded we have
\eq \limsup\limits_{t\ra\infty}\sup_{M}\frac{\vt'}{\vt H}\leq\frac{1}{n}.\eeq
\eLem

\pf
If $\vt'\leq c<\infty,$ this follows from Lemma \ref{HBound}. Thus suppose $\vt'\ra\infty.$ Define
\eq w=\log\left(\frac{1}{H}\right)+\log v+\log\vt'-\log\vt-\log\frac{1}{n}\eeq
and
\eq z=wf(u),\eeq
where $f\geq0$ will be determined later.
We have, as soon as $\vt'\geq 1,$

\begin{align}\begin{split}\label{HGrowth1} \dot{w}-\frac{1}{H^{2}}\Delta w&=-\frac{n}{H^{2}}\frac{\vt''}{\vt}-\frac{\t}{H^{2}}\|Du\|^{2}+\frac{1}{H^{2}}\frac{\Phi_{i}}{\Phi}\frac{\Phi^{i}}{\Phi}\\
				&\hp{=}+\frac{2}{H}\frac{\vt'}{\vt}v^{-1}-\frac{n}{H^{2}}\frac{\vt'^{2}}{\vt^{2}}+\frac{1}{\vt^{2}H^{2}}\|Du\|^{2}-\frac{1}{H^{2}}\frac{v_{i}}{v}\frac{v^{i}}{v}\\
				&\hp{=}+\frac{2}{H^{2}}\frac{\vt'}{\vt}\frac{v_{i}}{v}u^{i}+\left(\frac{\vt''}{\vt'}-\frac{\vt'}{\vt}\right)\frac{2}{H}v^{-1}-\left(\frac{\vt''}{\vt}-\frac{\vt'^{2}}{\vt^{2}}\right)\frac{n}{H^{2}}\\
				&\hp{=}+\left(\frac{\vt''}{\vt}-\frac{\vt'^{2}}{\vt^{2}}\right)\frac{1}{H^{2}}\|Du\|^{2}-\left(\frac{\vt''}{\vt'}-\frac{\vt'}{\vt}\right)'\frac{1}{H^{2}}\|Du\|^{2}\\
				&=-\frac{2n}{H^{2}}\frac{\vt''}{\vt}+\frac{2}{H}\frac{\vt''}{\vt'}v^{-1}+\frac{1}{H^{2}}\frac{\Phi_{i}}{\Phi}\frac{\Phi^{i}}{\Phi}-\frac{1}{H^{2}}\frac{v_{i}}{v}\frac{v^{i}}{v}\\
				&\hp{=}+\frac{2}{H^{2}}\frac{\vt'}{\vt}\frac{v_{i}}{v}u^{i}-\left(\frac{\vt''}{\vt'}-\frac{\vt'}{\vt}\right)'\frac{1}{H^{2}}\|Du\|^{2}.\end{split}\end{align}
Define \eq\~{w}=\sup\limits_{\xi\in M}w(\cdot,\xi)=w(t,\xi_{t}).\eeq
At the points $(t,\xi_{t})$ we have
\eq \frac{\Phi_{i}}{\Phi}=-\frac{v_{i}}{v}-\left(\frac{\vt''}{\vt'}-\frac{\vt'}{\vt}\right)u_{i}=vh^{k}_{i}u_{k}-\frac{\vt''}{\vt'}u_{i},\eeq
cf. \cite[(5.29)]{Gerhardt:01/1996}.
From (\ref{HGrowth1}) we obtain for almost every large $t$

\eq \dot{\~{w}}\leq \frac{2}{H^{2}}\frac{\vt''}{\vt'}\left(Hv^{-1}-n\frac{\vt'}{\vt}+c\frac{\vt'}{\vt}\|Du\|^{2}\right)\leq \frac{2}{H^{2}}\frac{\vt''}{\vt'}\left(Hv^{-1}-\frac{n}{2}\frac{\vt'}{\vt}\right),\eeq
using that $\|Du\|^{2}\ra 0$ in any case. Thus for a.e. large $t$ we have
\eq \dot{\~{w}}\leq \frac{2}{H}\frac{\vt''}{\vt'}v^{-1}\left(1-\frac{1}{2}e^{\~{w}}\right)\eeq
 and $w$ is bounded. Now let
 
 \eq f(u)=\vt'^{\g}, \g>0.\eeq
 Then from (\ref{HGrowth1}) we obtain
 
 \begin{align}\begin{split} \dot{z}-\frac{1}{H^{2}}\Delta z&\leq -\frac{2n}{H^{2}}\frac{\vt''}{\vt}\vt'^{\g}+\frac{2}{H}\frac{\vt''}{\vt'}v^{-1}\vt'^{\g}+\frac{1}{H^{2}}\frac{\Phi_{i}}{\Phi}\frac{\Phi^{i}}{\Phi}\vt'^{\g}-\frac{1}{H^{2}}\frac{v_{i}}{v}\frac{v^{i}}{v}\vt'^{\g}\\
 					&\hp{=}+\frac{2}{H^{2}}\frac{\vt'}{\vt}\frac{v_{i}}{v}u^{i}\vt'^{\g}-\left(\frac{\vt''}{\vt}-\frac{\vt'}{\vt}\right)'\frac{1}{H^{2}}\|Du\|^{2}\vt'^{\g}\\
					&\hp{=}+\frac{2\g}{H}\vt'^{\g-1}\vt''v^{-1}w-\frac{\g n}{H^{2}}\vt'^{\g}\frac{\vt''}{\vt}w\\
						&\hp{=}+\frac{\g}{H^{2}}\vt'^{\g}\frac{\vt''}{\vt}\|Du\|^{2}w-\g(\g-1)\vt'^{\g-2}\vt''^{2}\frac{1}{H^{2}}\|Du\|^{2}w\\
						&\hp{=}-\g\vt'^{\g-1}\vt'''\frac{1}{H^{2}}\|Du\|^{2}w-\frac{2}{H^{2}}w_{i}f^{i}.\end{split}\end{align}
						
Thus for \eq \~{z}=\sup_{\xi\in M}z(\cdot,\xi)\eeq
we obtain for almost every $t,$ also using
\eq wf_{i}=-w_{i}f=-\frac{\Phi_{i}}{\Phi}f-\frac{v_{i}}{v}f-\left(\frac{\vt''}{\vt'}-\frac{\vt'}{\vt}\right)u_{i}f,\eeq

\begin{align}\begin{split} \dot{\~{z}}&\leq \frac{2}{H^{2}}\frac{\vt''}{\vt'}f\left(Hv^{-1}-n\frac{\vt'}{\vt}+\g Hv^{-1}w\right)\\
					&\hp{=}+\frac{1}{H^{2}}\frac{\vt''}{\vt}\left(c\|Du\|^{2}\vt'^{\g}-\g n\~{z}+c(\g^{2}+\g)\|Du\|^{2}\vt'^{\g}\right)\\
					&=\frac{2n}{H^{2}}f\frac{\vt''}{\vt}\left(e^{-w}-1+\g we^{-w}\right)\\
					&\hp{=}+\frac{1}{H^{2}}\frac{\vt''}{\vt}\left(c\|Du\|^{2}\vt'^{\g}-\g n\~{z}+c\g(\g+1)\|Du\|^{2}\vt'^{\g}\).\end{split}\end{align}
Using Lemma \ref{GradDecay} and Lemma \ref{HDecay}  we see, that for small $\g>0,$ both brackets are negative on the set
\eq \{\~{z}\geq 1\}.\eeq
Thus $\~{z}$ is bounded and we conclude, that in case of unbounded $\vt'$
\eq \limsup\limits_{t\ra\infty}\sup_{M}\frac{\vt'}{\vt H}\leq \frac{1}{n}.\eeq
\epf

\subsection*{Optimal decay estimates}

Now we prove the decay estimates as formulated in Theorem \ref{mainthm}.

\Prop\label{optGrad}
Under the assumptions \ref{assA} and \ref{assB} let $u$ be the solution of (\ref{ScalarFlow}). Then there exist $c=c(N,M_{0})$ and $\mu>0,$ such that
\eq |Du|^{2}=v^{2}-1\leq ce^{-\mu t}\eeq and
\eq |Du|^{2}\leq \frac{c}{\vt'^{2}}\eeq
Furthermore, under the condition (\ref{RescaledtoSphere}) we have
\eq \vt'^{2}|Du|^{2}\ra 0.\eeq
\eProp

\pf
If $\vt'$ is uniformly bounded, this result was already proven in Lemma \ref{GradDecay} and Corollary \ref{GradDecayB}. Thus suppose $\vt'\ra\infty.$ We use both representations of the evolution equation for $v,$ namely (\ref{Evv}), which uses coordinates $(\xi^{i})\in M$ and (\ref{EvGrad1}), which uses coordinates $(x^{i})\in\S^{n}.$ They are connected via the time dependent diffeomorphisms
\eq (x^{i})=(x^{i}(t,\xi)).\eeq
In the notation of Lemma \ref{EvGrad} we have
\eq w=\frac{1}{2}|D\p|^{2}=\frac{1}{2}(v^{2}-1).\eeq
Using (\ref{Evv}), we obtain

\begin{align}\begin{split}\label{optGrad1} \dot{w}-\frac{1}{H^{2}}\Delta w&\leq -\frac{4}{H}\frac{\vt'}{\vt}\frac{v}{v+1}w-\frac{2n\t}{H^{2}}w+\frac{2}{\vt^{2}H^{2}}w\\
			&\hp{=}+\frac{c}{H^{2}}\|Dv\|^{2}+\frac{c}{H^{2}}\frac{\vt'}{\vt}\|Dv\|\|Du\|\\
			&\leq \(o(1)-\frac{2}{n}-\frac{2n}{H^{2}}\(\frac{\vt''}{\vt}-\frac{\vt'^{2}}{\vt^{2}}\)-\frac{2(n-1)}{\vt^{2}H^{2}}\)w\\
			&\hp{=}+\frac{c}{H^{2}}\|Dv\|^{2}+\frac{c}{H^{2}}\frac{\vt'}{\vt}\|Dv\|\|Du\|\\
			&\leq \(o(1)-\frac{2n}{H^{2}}\frac{\vt''}{\vt}-\frac{2(n-1)}{\vt^{2}H^{2}}\)w\\
			&\hp{=}+\frac{c}{H^{2}}\|Dv\|^{2}+\frac{c}{H^{2}}\frac{\vt'}{\vt}\|Dv\|\|Du\|\\
			&\leq2(n-1)\frac{\vt'^{2}}{\vt^{2}H^{2}}\(o(1)-\frac{\vt''\vt+1}{\vt'^{2}}\)w\\
			&\hp{=}+\frac{c}{H^{2}}\|Dv\|^{2}+\frac{c}{H^{2}}\frac{\vt'}{\vt}\|Dv\|\|Du\|  \end{split}\end{align}
Using \eqref{assB} we obtain the result on the exponential decay. To prove the second claim, come back to the proof of Lemma \ref{GradDecay} and consider (\ref{GradDecay3}) with $\g=2.$ We obtain for $z=w\vt'^{2},$ that for large $t$
  \eq \mc{L}z\leq ce^{-\mu t}z-\frac{c}{\vt'^{2}}z,\eeq
   where $c=c(\d,N,M_{0})$ and $\d$ is chosen so small, that we may absorb the derivative term in (\ref{GradDecay3}). However, $\d=\d(N,M_{0}).$ Thus $z$ is bounded. In case \eq \vt'^{-2}(\vt^{-1}(\vt(R_{0})e^{\frac{t}{n}}))\notin L^{1}([0,\infty)),\eeq we use Lemma \ref{BarriersB} and integration to show, that $\log z$ converges to $-\infty,$ which means that $z\ra 0.$ Thus the proof is complete.  
\epf

Now we show that the hypersurfaces become umbilic. We do this by estimating the largest principal curvature from above and $H$ from below.

\Lem\label{optHGrowth}
Under the assumptions \ref{assA} and \ref{assB} let $u$ be the solution of (\ref{ScalarFlow}). Then there exists $c=c(N,M_{0}),$ such that for all $t\in[0,\infty)$
\eq \frac{v}{H}\frac{\vt'}{\vt}-\frac{1}{n}\leq \frac{ct}{\vt'^{2}}\eeq and for all $\g<2$ there exists $c=c(\g,N,M_{0}),$ such that
\eq \frac{v}{H}\frac{\vt'}{\vt}-\frac{1}{n}\leq \frac{c}{\vt'^{\g}}\ \ \forall t\in[0,\infty).\eeq
\eLem

\pf
Using (\ref{HGrowth1}), we calculate the evolution equation for 
\eq z=\frac{v}{H}\frac{\vt'}{\vt},\eeq
namely

\begin{align}\begin{split} \dot{z}-\frac{1}{H^{2}}\Delta z&=-\frac{2n}{H^{2}}\frac{\vt''}{\vt}z+\frac{2}{H}\frac{\vt''}{\vt'}v^{-1}z+\frac{1}{H^{2}}\frac{\Phi_{i}}{\Phi}\frac{\Phi^{i}}{\Phi}z-\frac{1}{H^{2}}\frac{v_{i}}{v}\frac{v^{i}}{v}z\\
				&\hp{=}+\frac{2}{H^{2}}\frac{\vt'}{\vt}\frac{v_{i}}{v}u^{i}z-\left(\frac{\vt''}{\vt'}-\frac{\vt'}{\vt}\right)'\frac{1}{H^{2}}\|Du\|^{2}z-\frac{1}{H^{2}}w_{i}w^{i}z\\
						&\leq -\frac{2n}{H}v^{-1}\frac{\vt''}{\vt'}\(z-\frac{1}{n}\)z-\frac{2}{H^{2}}\frac{\Phi_{i}}{\Phi}\frac{v^{i}}{v}z\\
						&\hp{=}-\frac{2}{H^{2}}\(\frac{\vt''}{\vt'}-\frac{\vt'}{\vt}\)\frac{\Phi_{i}}{\Phi}u^{i}z+c\|Du\|^{2}z.\end{split}\end{align}
Define \eq \rho=\(z-\frac{1}{n}\)\vt'^{\g}\eeq
 and
 \eq \~{\rho}(t)=\sup_{M}\rho(t,\cdot)=\rho(t,\xi_{t}).\eeq
 We obtain for almost every $t,$ that
 \begin{align}\begin{split} \dot{\~{\rho}}&\leq -\frac{2n}{H}v^{-1}\frac{\vt''}{\vt'}z\~{\rho}-\frac{2}{H^{2}}\frac{\Phi_{i}}{\Phi}\frac{v^{i}}{v}z\vt'^{\g}-\frac{2}{H^{2}}\(\frac{\vt''}{\vt'}-\frac{\vt'}{\vt}\)\frac{\Phi_{i}}{\Phi}u^{i}z\vt'^{\g}\\
 				&\hp{=}+c\|Du\|^{2}\vt'^{\g}+\frac{2\g}{H}\frac{\vt''}{\vt'}v^{-1}\~{\rho}-\frac{\g n}{H^{2}}\frac{\vt''}{\vt}\~{\rho}-\frac{2}{H^{2}}z_{i}(\vt'^{\g})^{i}\\
				&=-\frac{2n}{H}v^{-1}\frac{\vt''}{\vt'}\~{\rho}\(z+\frac{\g}{2}z-\frac{\g}{n}\)+c\|Du\|^{2}\vt'^{\g},\end{split}\end{align}
where we also used
\eq \rho_{i}=0\eeq
at $(t,\xi_{t}).$
Thus Proposition \ref{optGrad} yields in case $\g=2,$ that
\eq \~{\rho}(t)\leq ct\ \ \forall t\in[0,\infty)\eeq
and in case $\g<2$
\eq \~{\rho}\leq c.\eeq
\epf

\Lem\label{kappaDecay}
Under the assumptions \ref{assA} and \ref{assB} let $u$ be the solution of (\ref{ScalarFlow}). Then there exists $c=c(N,M_{0}),$ such that the principal curvatures of the flow hypersurfaces satisfy
\eq \k_{i}\frac{\vt}{\vt'}\leq c\eeq
 and in case $\vt'\ra\infty$ we have
 \eq \left|v\k_{i}\frac{\vt}{\vt'}-1\right|\leq \frac{ct}{\vt'^{2}}.\eeq
 Furthermore, for all $\g<2$ there exists $c=c(\g,N,M_{0}),$ such that 
 \eq \left|v\k_{i}\frac{\vt}{\vt'}-1\right|\leq \frac{c}{\vt'^{\g}}.\eeq
\eLem

\pf
The first part was already shown in Lemma \ref{KappaBound}. Thus suppose $\vt'\ra\infty.$ Consider (\ref{KappaBound1}) with 
\eq f(v)=\log v,\eeq
such that
\eq z=\log\zeta+\log\frac{\vt}{\vt'}+\log v.\eeq
We want to bound the function 
\eq \rho=g(t)(e^{z}-1)\vt'^{\g},\eeq
where $g$ and $\g$ are to be determined. Thus suppose
\eq 1\leq\sup_{[0,T]\times M}\rho=\rho(t_{0},\xi_{0}), t_{0}>0.\eeq
Again, without loss of generality assume
\eq g_{ij}=\d_{ij}, h_{ij}=\k_{i}\d_{ij}, \k_{1}\leq \dots\k_{n}.\eeq
(\ref{KappaBound1}) becomes
\begin{align}\begin{split} \dot{z}-\frac{1}{H^{2}}\Delta z&\leq -\frac{2}{H}\(\k_{n}-\frac{\vt'}{\vt}\)+c\|Du\|^{2}+\frac{1}{H^{2}}(\log h^{n}_{n})_{i}(\log h^{n}_{n})^{i}\\
		&\hp{=}+\frac{2}{H}v^{-1}\(\frac{\vt'}{\vt}-\frac{\vt''}{\vt'}\)\(1-v^{-1}\k_{n}^{-1}\frac{\vt'}{\vt}\).\end{split}\end{align}
Thus 
\eq w=e^{z}\eeq
satisfies
\begin{align}\begin{split} \dot{w}-\frac{1}{H^{2}}\Delta w&\leq -\frac{2}{H}\frac{\vt'}{\vt}(w-1)^{2}-\frac{2}{H}\frac{\vt''}{\vt'}v^{-1}(w-1)+c\|Du\|^{2}w\\
			&\hp{=}+\frac{c}{H^{2}}\frac{\vt'}{\vt}|(\log h^{n}_{n})_{i}u^{i}|,\end{split}\end{align}
so that
\begin{align}\begin{split} \dot{\rho}-\frac{1}{H^{2}}\Delta\rho&\leq -\frac{2}{H}\frac{\vt'}{\vt}(w-1)\rho-\frac{2}{H}\frac{\vt''}{\vt'}v^{-1}\rho+cg\|Du\|^{2}\vt'^{\g}\\
			&\hp{=}+\frac{cg}{H^{2}}\frac{\vt'}{\vt}|(\log h^{n}_{n})_{i}u^{i}|\vt'^{\g}+\frac{2\g}{H}\frac{\vt''}{\vt'}v^{-1}\rho-\frac{\g n}{H^{2}}\frac{\vt''}{\vt}\rho\\
			&\hp{=}-\frac{2g}{H^{2}}w_{i}(\vt'^{\g})^{i}+\frac{g'}{g}\rho.\end{split}\end{align}
Using, that at $(t_{0},\xi_{0})$ there holds
\begin{align}\begin{split} 0=\frac{\rho_{i}}{g}&=w_{i}\vt'^{\g}+\g\vt'^{\g-1}\vt''(w-1)u_{i}\\
						&=h^{n}_{n;i}\vt\vt'^{\g-1}v+h^{n}_{n}\(\frac{\vt}{\vt'}\)_{i}v\vt'^{\g}\\
						&\hp{=}+h^{n}_{n}\frac{\vt}{\vt'}v_{i}\vt'^{\g}+\g\vt'^{\g-1}\vt''(w-1)u_{i},\end{split}\end{align}
we obtain
\eq \|D(\log h^{n}_{n})\|\leq c\frac{\vt'}{\vt}\|Du\|.\eeq
Thus
\begin{align}\begin{split} \label{kappaDecay1}0&\leq -\frac{2}{H}\frac{\vt'}{\vt}(w-1)\rho-\frac{\g n}{H}\frac{\vt''}{\vt'}v^{-1}\rho\(\frac{v}{H}\frac{\vt'}{\vt}-\frac{2\g-2}{\g n}\)\\
			&\hp{=}+cg\|Du\|^{2-\g}+\frac{g'}{g}\rho.\end{split}\end{align}
Let $\mu>0$ be as in Proposition \ref{optGrad}. First, let $0<\g<2$ and
\eq g(t)=\(\int_{-1}^{t}e^{-\e s}ds\)^{-1},\  0<\e<\frac{1}{2}(2-\g)\mu.\eeq
Then from Lemma \ref{HDecay} we obtain for large $t_{0},$
\eq 0\leq cge^{-\frac{(2-\g)\mu}{2}t_{0}}-ge^{-\e t_{0}}\rho<0.\eeq
Since $g$ is decreasing but still strictly positive, we find
\eq \label{kappaDecay2} (w-1)\vt'^{\g}\leq c(\g,N,M_{0}),\eeq
as well as
\eq \label{kappaDecay3}\frac{H}{v}\frac{\vt}{\vt'}=\sum_{i=1}^{n}\(v\k_{i}\frac{\vt}{\vt'}-1\)+n\leq c\vt'^{-\g}+n.\eeq
Consider (\ref{kappaDecay1}) with $\g=2$ and $g(t)=t^{-1}.$
Then
\begin{align}\begin{split} 0&\leq \frac{2}{H^{2}}\frac{\vt''}{\vt}\(\frac{H}{v}\frac{\vt}{\vt'}-n\)\rho+\frac{1}{t}(c-\rho)\\
			&=\frac{2}{tH^{2}}\frac{\vt''}{\vt}(w-1)\(\frac{H}{v}\frac{\vt}{\vt'}-n\)\vt'^{2}+\frac{1}{t}(c-\rho)\\
			&\leq \frac{2c}{tH^{2}}\frac{\vt''}{\vt}(w-1)\vt'+\frac{1}{t}(c-\rho)\\
			&<0\end{split}\end{align}
for large $\rho.$ Here we used (\ref{kappaDecay2}) and (\ref{kappaDecay3}).
 From Lemma \ref{optHGrowth} we get
\begin{align}\begin{split} \frac{\sum_{i=1}^{n}\(1-v\k_{i}\frac{\vt}{\vt'}\)}{nvH\frac{\vt}{\vt'}}&=\frac{n-vH\frac{\vt}{\vt'}}{nvH\frac{\vt}{\vt'}}\\
		&\leq \frac{ct}{\vt'^{2}}\ \ \(\frac{c}{\vt'^{\g}}\ \text{respectively}\)\end{split}\end{align}
and thus

\begin{align}\begin{split} 1-v\k_{1}\frac{\vt}{\vt'}&\leq \frac{ct}{\vt'^{2}}+c\sum_{i=2}^{n}\(v\k_{i}\frac{\vt}{\vt'}-1\)\\
					&\leq\frac{ct}{\vt'^{2}}\ \ \(\frac{c}{\vt'^{\g}}\ \text{respectively}\).\end{split}\end{align}
\epf

Hence, we have proven part (ii) of Theorem \ref{mainthm} completely. We will finish this paper by proving the optimal decay of the second fundamental form in case of improving pinching of $N.$

\begin{thm}\label{optKappa}
Under the assumptions \ref{assA}, \ref{assB} and \ref{assC} let $u$ be a solution of (\ref{ScalarFlow}). Then there exists $c=c(N,M_{0}),$ such that
\eq \left|h^{i}_{j}-\frac{\vt'}{\vt}\d^{i}_{j}\right|\leq \frac{c}{\vt'}e^{-\frac{t}{n}}.\eeq
\end{thm}

\pf
Only the case $\vt'\ra\infty$ has to be considered. From (\ref{EvA}) and (\ref{Evu}) we obtain, that the function
\eq G=\frac{1}{2}\left\|A-\frac{\vt'}{\vt}I\right\|^{2}=\frac{1}{2}\(h^{i}_{j}-\frac{\vt'}{\vt}\d^{i}_{j}\)\(h^{j}_{i}-\frac{\vt'}{\vt}\d^{j}_{i}\)\eeq
satisfies

\begin{align}\begin{split}	&\dot{G}-\frac{1}{H^{2}}\Delta G\\
				=&\frac{1}{H^{2}}\(\|A\|^{2}-2\frac{\vt'}{\vt}H+n\frac{\vt'^{2}}{\vt^{2}}\)\(\|A\|^{2}-\frac{\vt'}{\vt}H\)\\
				&-\frac{2}{H}\mrm{tr}\(A-\frac{\vt'}{\vt}I\)^{3}-\frac{n}{\vt^{2}H^{2}}\(\|A\|^{2}-\frac{\vt'}{\vt}H\)-\frac{4}{H}\frac{\vt'}{\vt}G\\
				&+\frac{(n+1)\t}{H^{2}}\|Du\|^{2}\(\|A\|^{2}-\frac{\vt'}{\vt}H\)-\frac{2n\t}{H^{2}}G\\
				&-\frac{n}{\vt^{2}H^{2}}\frac{\vt'}{\vt}\(H-n\frac{\vt'}{\vt}\)-\frac{2\t}{H}\frac{v-1}{v^{2}}\(H-n\frac{\vt'}{\vt}\)\\
				&+\frac{2}{\vt^{2}H}v^{-1}\(H-n\frac{\vt'}{\vt}\)+\frac{\t'}{vH^{2}}(nu_{j}u^{i}-\|Du\|^{2}\d^{i}_{j})h^{j}_{i}\\
				&+\frac{n\t}{H^{2}}\(h^{i}_{m}u^{m}u_{j}+h^{m}_{j}u_{m}u^{i}-\frac{2}{n}h^{m}_{k}u_{m}u^{k}\d^{i}_{j}\)h^{j}_{i}\\
				&+\frac{2\t}{H^{2}}v^{-1}\frac{\vt'}{\vt}(\|Du\|^{2}\d^{i}_{j}-nu^{i}u_{j})h^{j}_{i}\\
				&-\frac{2}{H^{3}}H_{j}H^{i}\(h^{j}_{i}-\frac{\vt'}{\vt}\d^{j}_{i}\)-\(\frac{\vt''}{\vt}-\frac{\vt'^{2}}{\vt^{2}}\)\frac{\vt'}{\vt}\frac{\|Du\|^{2}}{H^{2}}\(H-n\frac{\vt'}{\vt}\)\\
				&+\(\frac{\vt'}{\vt}\)''\frac{\|Du\|^{2}}{H^{2}}\(H-n\frac{\vt'}{\vt}\)-\frac{1}{H^{2}}\left\|D\(A-\frac{\vt'}{\vt}I\)\right\|^{2}.\end{split}\end{align}
Now define \eq w=Gf(u),\eeq
where $f>0$ will be defined later. Then, as soon as $\vt'\geq 1,$

\begin{align}\begin{split} \label{optKappa1}	&\dot{w}-\frac{1}{H^{2}}\Delta w\\
					\leq& -\frac{2}{H}\mrm{tr}\(A-\frac{\vt'}{\vt}I\)^{3}f+\frac{2}{H^{2}}h^{i}_{j}\(h^{j}_{i}-\frac{\vt'}{\vt}\d^{j}_{i}\)w-\frac{4}{H}\frac{\vt'}{\vt}w\\
						&-\frac{n}{\vt^{2}H^{2}}\(h^{i}_{j}+\frac{\vt'}{\vt}\d^{i}_{j}-\frac{2}{nv}H\d^{i}_{j}\)\(h^{j}_{i}-\frac{\vt'}{\vt}\d^{j}_{i}\)f-\frac{2n\t}{H^{2}}w\\
						&+\frac{c|\t|}{H^{2}}\|Du\|^{2}\frac{\vt'}{\vt}\left\|A-\frac{\vt'}{\vt}I\right\|f+\frac{c|\t'|}{H^{2}}\|Du\|^{2}\left\|A-\frac{\vt'}{\vt}I\right\|f\\
						&+\frac{c\vt'}{\vt^{3}H^{2}}\|Du\|^{2}\|A-\frac{\vt'}{\vt}\mrm{I}\|f\\
						&-\frac{2}{H^{3}}H_{j}H^{i}\(h^{j}_{i}-\frac{\vt'}{\vt}\d^{j}_{i}\)f+\frac{2}{H}\frac{f'}{f}v^{-1}w-\frac{n}{H^{2}}\frac{\vt'}{\vt}\frac{f'}{f}w\\
						&+\frac{1}{H^{2}}\frac{\vt'}{\vt}\frac{f'}{f}\|Du\|^{2}w-\frac{1}{H^{2}}\frac{f''}{f}\|Du\|^{2}w-\frac{2}{H^{2}}\frac{G_{k}}{G}\frac{f^{k}}{f}w\\
						&-\frac{1}{\vt'^{2-\d}H^{2}}\(h^{i}_{j}-\frac{\vt'}{\vt}\d^{i}_{j}\)_{;k}{\(h^{j}_{i}-\frac{\vt'}{\vt}\d^{j}_{i}\)_{;}}^{k}f,\end{split}\end{align}
$0<\d<2$ to be chosen appropriately.
There holds
\begin{align}\begin{split}	\left\|D\(A-\frac{\vt'}{\vt}I\)\right\|^{2}&\geq \|DA\|^{2}-2\left|\(\frac{\vt'}{\vt}\)'\right|\|DH\|\|Du\|\\
					&\hp{=}-n\(\frac{\vt'}{\vt}\)'^{2}\|Du\|^{2}\\
					&\geq \|DA\|^{2}-\e\left|\(\frac{\vt'}{\vt}\)'\right|\frac{\vt^{2}}{\vt'^{2}}\|DH\|^{2}\\
					&\hp{=}-\frac{1}{\e}\left|\(\frac{\vt'}{\vt}\)'\right|\frac{\vt'^{2}}{\vt^{2}}\|Du\|^{2}-n\left|\(\frac{\vt'}{\vt}\)'\right|^{2}\|Du\|^{2}\\
					&\geq \frac{1}{2n}\|DH\|^{2}-\frac{c}{\e}\frac{\vt'^{4}}{\vt^{4}}\|Du\|^{2},\ 0<\e<<1.\end{split}\end{align}
Furthermore, from Proposition \ref{optGrad} and Lemma \ref{kappaDecay}, we obtain
\eq \left\|A-\frac{\vt'}{\vt}I\right\|\leq \frac{c}{\vt'^{\a}}\frac{\vt'}{\vt},\ \a<2.\eeq
Setting
\eq f=\vt^{\g},\ 2<\g<\min\left\{\l,2+\frac{\mu}{2}\right\},\eeq
where $\l$ is as in Assumption \ref{assC}, $\d=1$ and $\a>1,$ we obtain from (\ref{optKappa1}) at a point $(t_{0},\xi_{0})$ with
$ \sup_{[0,T]\times M}w=w(t_{0},\xi_{0})\geq 1,$
$t_{0}>0$ large, that

\begin{align}\begin{split} 0&\leq w\(o(1)-\frac{4}{H}\frac{\vt'}{\vt}+c\vt^{1-\l+\frac{\g}{2}}\vt'^{-1}+\frac{2\g}{H}\frac{\vt'}{\vt}v^{-1}-\frac{\g n}{H^{2}}\frac{\vt'^{2}}{\vt^{2}}\)\\
			&\hp{=}-\frac{n}{\vt^{2}H^{2}}\(h^{i}_{j}+\frac{\vt'}{\vt}\d^{i}_{j}-\frac{2}{nv}H\d^{i}_{j}\)\(h^{j}_{i}-\frac{\vt'}{\vt}\d^{j}_{i}\)\vt^{\g}\\
			&\hp{=}+\frac{c}{H^{2}}\frac{\vt^{\g}}{\vt'^{\a}}\|DH\|^{2}-\frac{1}{\vt'H^{2}}\frac{1}{2n}\|DH\|^{2}\vt^{\g}+\frac{c}{\e H^{2}}\frac{\vt'^{3}}{\vt^{4}}\vt^{\g}\|Du\|^{2}\\
			&\leq \(o(1)+(\g-4)\frac{1}{H}\frac{\vt'}{\vt}\)w+\frac{2}{\vt^{2}H^{2}}v^{-1}\(H-n\frac{\vt'}{\vt}\)^{2}\vt^{\g}\\
			&\hp{=}+\frac{cn}{\vt^{2}H^{2}}\|Du\|^{2}\frac{\vt'}{\vt}\left|H-n\frac{\vt'}{\vt}\right|\vt^{\g}+\frac{c}{H^{2}}\frac{\vt'^{\g}}{\vt'^{\a}}\|DH\|^{2}\\
			&\hp{=}-\frac{1}{\vt'H^{2}}\frac{1}{2n}\|DH\|^{2}\vt^{\g}+\frac{c}{\e H^{2}}\frac{\vt'^{2}}{\vt^{2}}e^{-\frac{\mu}{2}t_{0}}e^{(\g-2)\frac{t_{0}}{n}}\\
			&\leq \(o(1)+(\g-4)\frac{1}{H}\frac{\vt'}{\vt}\)w+ce^{\(\g-2-\frac{\mu}{2}\)t_{0}}<0. \end{split}\end{align}
Thus we have
\eq \left\|A-\frac{\vt'}{\vt}I\right\|\leq \frac{c}{\vt^{\frac{\g}{2}}}=\frac{c}{\vt'\vt^{\frac{\g}{2}-1}}\frac{\vt'}{\vt}.\eeq
Now consider (\ref{optKappa1}) with $f=\vt^{2}\vt'^{2},$ $0<\d<\min\{\l-2,\mu\}$ and $2>\a>2-\frac{\d}{2}.$ Then the function $\~{w}=\sup_{\xi\in M}w(\cdot,\xi)$
satisfies for almost every large $t,$ that

\begin{align}\begin{split} \dot{\~{w}}&\leq \Big(ce^{\(1-\frac{\g}{2}\)\frac{t}{n}}-\frac{4}{H}\frac{\vt'}{\vt}+ce^{(2-\l)t}+\frac{4}{H}\frac{\vt'}{\vt}v^{-1}+\frac{4}{H}\frac{\vt''}{\vt'}v^{-1}\\
					&\hp{\leq w(c}-\frac{2n}{H^{2}}\frac{\vt'^{2}}{\vt^{2}}-\frac{2n}{H^{2}}\frac{\vt''}{\vt}+ce^{-\mu t}\Big)\~{w}+ce^{(2-\g)\frac{t}{n}}\\
					&\hp{=}+\frac{2}{\vt^{2}H^{2}}v^{-1}\(H-n\frac{\vt'}{\vt}\)^{2}\vt^{2}\vt'^{2}+\frac{cn}{\vt^{2}H^{2}}\|Du\|^{2}\frac{\vt'}{\vt}\left|H-n\frac{\vt'}{\vt}\right|\vt^{2}\vt'^{2}\\
					&\hp{=}+\frac{2}{H^{3}}\|DH\|^{2}\left\|A-\frac{\vt'}{\vt}I\right\|\vt'^{2}\vt^{2}-\frac{1}{H^{2}}\|DA\|^{2}\vt'^{\d}\vt^{2}\\
					&\hp{=}+\frac{2|\t|}{H^{2}}\|DH\|\|Du\|\vt^{2}\vt'^{\d}+\frac{2}{\vt^{2}H^{2}}\|DH\|\|Du\|\vt'^{\d}\vt^{2}\\
					&\hp{=}+\frac{n\t^{2}}{H^{2}}\|Du\|^{2}\vt'^{\d}\vt^{2}+\frac{2n|\t|}{\vt^{2}H^{2}}\|Du\|^{2}\vt'^{\d}\vt^{2}+\frac{n}{\vt^{4}H^{2}}\|Du\|^{2}\vt'^{\d}\vt^{2}\\
					&\leq ce^{-\e t}(w+1)+\frac{\e|\t|}{H^{2}}\frac{\vt^{\l}}{\vt'^{2}}\|DH\|^{2}\vt'^{\d}\vt^{2}+\frac{|\t|}{\e H^{2}}\frac{\vt'^{2}}{\vt^{\l}}\|Du\|^{2}\vt'^{\d}\vt^{2}\\
					&\hp{=}+\frac{2}{\vt^{\d}H^{2}}\|DH\|^{2}\vt'^{\d}\vt^{2}+\frac{2}{\vt^{2}H^{2}}\vt^{\d-2}\|Du\|^{2}\vt'^{\d}\vt^{2}\\
					&\hp{=}+\frac{cn}{H^{2}}\frac{\vt'^{4}}{\vt^{2\l}}\|Du\|^{2}\vt'^{\d}\vt^{2}+\frac{c}{H^{2}}\|DH\|^{2}\vt'^{\frac{\d}{2}}\vt^{2}+c\frac{\vt'^{2}}{\vt^{\l}}\vt^{2}\vt'^{\d-2}\|Du\|^{2}\\
					&\hp{=}+c\|Du\|^{2}\vt'^{\d-2}-\frac{1}{H^{2}}\|DA\|^{2}\vt'^{\d}\vt^{2}\\
					&\leq ce^{-\e t}(\~{w}+1),\end{split}\end{align}
	for small $\e>0$ and large $t.$
Thus $w$ is bounded and the proof complete.
\epf

\section{Concluding remarks}
We have seen, that under rotational symmetry of the ambient manifold, the IMCF does not produce singularities, if the initial hypersurface is starshaped. We should mention that one can say significantly more, if we rule out Euclidean behavior in the sense, that the sectional curvatures vary between strictly negative numbers. Then we are optimistic, that Assumption \ref{assB} may be weakened, namely (\ref{assB2}) removed.

The regularity assumptions on $N$ and on the initial hypersurface may surely be weakened. For an overview over weaker sufficient regularity assumptions have a look at \cite[Theorem 2.5.19]{Gerhardt:/2006}. Since the methods applied in the proofs of our work do not differ in case of weaker initial regularity, we did not consider it to be worthy formulating everything in the generality found in \cite[Theorem 2.5.19]{Gerhardt:/2006}.
Besides that, it should be possible to generalize the results of this work to other curvature functions of homogeneity 1. The higher homogeneity case, e.g. the Gaussian curvature, causes another difficulty, since even in the Euclidean setting we do not have long time existence. Here, other methods had to be applied, compare \cite{Gerhardt:01/2014}, those of which we are not quite aware, whether they are applicable in our general setting.

Furthermore, we would find it interesting, whether it would be possible to remove the extra $t$ factor in (\ref{umbilic}) even without Assumption \ref{assC}. It seems to us, that it might not, because the extra terms in the evolution equations of the curvature quantities do in general decay as fast as $\|Du\|^{2},$ so that in the rescaled evolution equations they are of order zero and do not decay. However, we could not come up with a counterexample. 

Although we are optimistic, that the results at hand should be enough for many applications regarding geometric inequalities, it would still be interesting to derive higher oder estimates in order to show that the rescaled surfaces $\~{u}=\log \vt-\frac{t}{n}$ are bounded in $C^{\infty}.$

\bibliographystyle{hamsplain}
\bibliography{Bibliography}

\end{document}